\documentclass{amsart}
\usepackage{graphicx}
\usepackage{psfrag}

\newtheorem{thm}{Theorem}[section]
\newtheorem{cor}[thm]{Corollary}
\newtheorem{lem}[thm]{Lemma}

\newtheorem{prop}[thm]{Proposition}

\theoremstyle{definition}
\newtheorem{defn}[thm]{Definition}
\newtheorem{rem}[thm]{Remark}

\numberwithin{equation}{section}

\newcommand{\Ss}{\mathbb S}
\newcommand{\Zz}{\mathbb Z}
\newcommand{\Nn}{\mathbb N}

\newcommand{\ie}{{\it i.e.}}

\newcommand{\surf}{{\mathcal F}} 
\newcommand{\sinfty}{\Ss^1_{\infty}} 

\newcommand{\aaxis}{{\tt A}}
\newcommand{\axis}[1]{\aaxis_{#1}}
\def\ug{\ucover\gamma}

\def\gg{{\gamma}}  

\def\rtp{r}  
\def\ltp{l}  

\def\Re{\mbox{Re}}
\def\Im{\mbox{Im}}

\def\Cc{{\mathbb C}} 

\def\ucover#1{\widetilde{#1}} 

\def\ug{\ucover\gamma}

\def\Rr{{\mathbb R}} 

\def\ddefine#1{{\em #1}}

\let\mgp=\marginpar \def\marginpar#1{\mgp{\raggedright\tiny #1}}
\def \cmp#1{\marginpar{\%\% \ #1 ---c}}

\let\lbl=\label
\def\label#1{\lbl{#1}\ifinner\else\marginpar{\ref{#1} #1}\ignorespaces\fi}

\def \marginpar#1{}

\begin{document}


\psfrag{p}{$p$}
\psfrag{q}{$q$}
\psfrag{x}{$x$}
\psfrag{r}{$r$}
\psfrag{y}{$y$}
\psfrag{b}{$\beta$}
\psfrag{a}{$\alpha$}
\psfrag{A}{$\aaxis$}
\psfrag{s}{$\sigma$}
\psfrag{C}{$C$}
\psfrag{c}{$c$}
\psfrag{l}{$l$}

\psfrag{b1}{$\beta_1$}
\psfrag{b2}{$\beta_2$}

\psfrag{Af}{$\axis{f}$}
\psfrag{finf}{$f(\infty)$}
\psfrag{fugx}{$f(\ug_x)$}
\psfrag{fugr}{$f(\ug_r)$}
\psfrag{f2ugr}{$f^2(\ug_r)$}
\psfrag{fugc}{$f(\ug_c)$}
\psfrag{fc}{$f(c)$}
\psfrag{finvugr}{$f^{-1}(\ug_r)$}
\psfrag{finv2ugr}{$f^{-2}(\ug_r)$}
\psfrag{fink}{$f^{-1}(\infty)$}
\psfrag{finkus}{$f^{-2}(\infty)$}
\psfrag{ginf}{$g(\infty)$}
\psfrag{ga}{$g(a)$}
\psfrag{gb1}{$g(\beta_1)$}
\psfrag{ginvp}{$g^{-1}(p)$}
\psfrag{gink}{$g^{-1}(c)$}
\psfrag{ginq}{$g^{-1}(q)$}
\psfrag{ginvugc}{$g^{-1}(\ug_c)$}
\psfrag{gp}{$g(p)$}
\psfrag{gr}{$\gamma_r$}
\psfrag{gx}{$g(x)$}
\psfrag{gugx}{$g(\ug_x)$}
\psfrag{gugc}{$g(\ug_c)$}
\psfrag{gugr}{$g(\ug_r)$}
\psfrag{gug0}{$g(\ug_0)$}
\psfrag{gug1}{$g(\ug_1)$}
\psfrag{h1ugfinf}{$h_1(\ug_{f(\infty)})$}
\psfrag{hink}{$h_2(\ug_x)$}

\psfrag{ldots}{$\ldots$}

\psfrag{x0ex}{$x_0=x$}
\psfrag{x1efx0}{$x_1=f(x_0)$}
\psfrag{x2}{$x_2$}
\psfrag{rrr}{$r=f(r)$}
\psfrag{ug1efug0}{$\ug_1=f(\ug_0)$}
\psfrag{ug0eugx}{$\ug_0=\ug_x$}

\psfrag{ug2}{$\ug_2$}
\psfrag{ugr}{$\ug_r$}
\psfrag{ugfinf}{$\ug_{f(\infty)}$}
\psfrag{ugx}{$\ug_x$}
\psfrag{ugc}{$\ug_c$}
\psfrag{uginf}{$\ug_{g(\infty)}$}


\title[Simple geodesics]{Simple geodesics on a punctured surface}%
\author{Chaim Goodman--Strauss and Yo'av Rieck}%
\address{Department of mathematical Sciences, University of
Arkansas, Fayetteville, AR 72701}
\email{strauss@uark.edu}%
\email{yoav@uark.edu}

\thanks{The first named author was supported in part by NSF through
grant number DMS-0072573  The second named author was supported in
part by the 21st century COE program ``Constitution for wide-angle
mathematical basis focused on knots"
(Osaka City University); leader: Akio Kawauchi.}%

\date{June 8 2005}%

\begin{abstract}
We give a new proof of McShane's Theorem \cite{mcshane}, using
simple equivariant methods  in the hyperbolic plane.
\end{abstract}
\maketitle

\section{Introduction and statements of results}
\label{sec:intro}

Let $\surf$ be a  orientable hyperbolic surface, complete and of
finite area, with one cusp (denoted $c$). Let $N(c)$ be a normal
neighborhood of $c$, and let $\Ss^1 = \partial N(c)$.  An oriented
geodesic is called \em cuspidal \em if it comes out of the cusp.
We study simple cuspidal oriented geodesics on $\surf$. Let $\Ss^1
= \partial N(c)$ be a circle that bounds the cusp, chosen
symmetrically (see next section).  Then we may parameterize the
oriented geodesics out of the cusp by $\Ss^1$: every oriented
cuspidal geodesic determines a unique point on $\Ss^1$ that is its
first intersection with $\Ss^1$.  For a point $x \in \Ss^1$ we
denote this geodesic by $\gg_x$.  We say that $\gamma_x$ is \em
bicuspidal \em if and only if it has both its ends in the cusp (in
that case there is a point $x \neq y \in \Ss^1$ so that $\gamma_y$
is $\gamma_x$ with reverse orientation); otherwise $\gamma_x$ is
called \em unicuspidal\em. Let $E \subset \Ss^1$ be all the points
$x$ for which $\gamma_x$ is simple. McShane proved the following
theorem (\cite[Theorem 4]{mcshane}):

\begin{thm}[McShane]
\label{thm:mcshane}

Let $\surf$ be a complete hyperbolic surface with finite area and
a cusp $c$. With notation as above, $E$ consists of a Cantor set
(say $K$) and isolated points, so that the following holds:

\begin{enumerate}
    \item $x \subset E$ is isolated if and only if $\gamma_x$
        is bicuspidal;
    \item $x \subset K$ is an endpoint of $K$
    if and only if $\gamma_x$ spirals onto a simple closed
    curve.\footnote{The complement of a  Cantor set
    $K\subset\Ss$ is a collection of open intervals;
    by ``the endpoints of $K$" we mean  the endpoints
    of these intervals.  $K$ has a countably infinite
    set of endpoints and uncountable set of points
    that are not endpoints.}
\end{enumerate}

Moreover, every connected component of $\Ss^1 \setminus K$
contains exactly one isolated point of $E$.
\end{thm}

Thus a component of $\Ss^1 \setminus K$ (say $J$) is bound by two
points (say $x$ and $y$) so that $\gg_x$ and $\gg_y$ spiral onto
simple closed geodesics, say $\alpha$ and $\beta$.  It is easy to
see that these geodesics co-bound a pair of pants with the cusp.
By using basic hyperbolic trigonometry \cite[Proposition
3]{mcshane} McShane showed that the length of $J$ is
${l(\Ss^1)}/({1+e^\frac{l(\alpha)+l(\beta)}{2}})$, where
$l(\cdot)$ denotes the hyperbolic length restricted to $\Ss^1$.
This, and the fact that by Birman and Series \cite{birman-series}
the measure of $E$ is zero, implies:

\begin{thm}[McShane's Identity]
\label{thm;mcshane's identity}

$$\Sigma_{\alpha,\beta}\frac{1}{1+e^\frac{|\alpha|+|\beta|}{2}} = 1/2$$

\noindent where the sum is taken over all pairs of geodesics
$\alpha$, $\beta$ that co-bound a pair of pants with the cusp $c$.
\end{thm}

Many authors worked on McShane's Identity (see references)
including Akiyoshi, Skakuma and Yamashita; Bowditch; Mirzakhani;
Tan, Wong and Zhang. Their work resulted in many alternate proofs
and generalizations 
of Theorem~\ref{thm;mcshane's identity}.

The purpose of this article is to give a short elementary proof of
Theorem~\ref{thm:mcshane}.  The proof is set in the universal
cover of $\surf$ (denoted $\ucover{\surf}$ and identified with the
upper half plane in $\Cc$) and uses basic geometry of the
hyperbolic plane and the action of $\pi_1(\surf)$ on
$\ucover{\surf}$; we view $\pi_1(\surf)$ as a subgroup of the
isometries of $\ucover{\surf}$.  Our proof is based on a new
interpretation of the open intervals complementary to $K$; McShane
refers to these intervals as \em gaps\em; below we call them \em
deadzones\em. We remark that a geodesic on $\surf$ is simple if
and only if any two of its lifts to $\ucover{\surf}$ are disjoint.
We call a geodesic in $\ucover{\surf}$ \em simple \em if it is
disjoint from all it images under $\pi_1(\surf)$; thus the
statements ``$\gg \subset \surf$ is simple" and ``some lift of
$\gg$ to $\ucover\surf$ is simple" and ``all lifts of $\gg$ to
$\ucover{\surf}$ are simple" are equivalent.

\noindent {\bf Remark.} It is well known that McShane's Theorem
holds for surfaces with more than one cusp and with totally
geodesic boundary components (see, for example, \cite{tan-main}).
Our techniques can be generalized to those settings as well.
However, as it is our goal to give a \em\ simple \em\ proof we
will not do so here.  

\noindent{\bf Acknowledgements:}  We would like to thank Max
Forester for teaching us about McShane's
work and the body of work related to it.  We also thank Makoto
Sakuma and Yasushi Yamashita for helpful conversations.

\section{Background: planar hyperbolic geometry}
\label{sec:background}

Though we assume some familiarity with the basic notions of
hyperbolic geometry, we pause to sketch out some essential
terminology, notations and ideas; the reader familiar with planar
hyperbolic geometry may skip this section as the results in the
section are well-known.  The hyperbolic plane may be compactified
by adding a circle at infinity that we denote $\sinfty$ (defined
by equivalence classes of geodesics that approach each other
asymptotically).  We work in the upper half plane model which we
identify with $\{z \in \Cc | \Re(z) > 0\}$. Then we get that
$\sinfty = \mathbb R \cup \{\infty\}$. Geodesics in
$\ucover{\surf}$ are in one-to-one correspondence with pairs of
distinct points in $\sinfty$ and oriented geodesics are in
one-to-one correspondence with ordered distinct pairs of points.

In the upper-half plane model, we arrange $\ucover{\surf}$ so that
$\infty$ is a lift of the cusp and its cuspidal subgroup is
generated by $z \to z+1$. Then every cuspidal oriented geodesic
has a vertical lift (clearly not unique) oriented downwards and
ending at some $x \in \Rr$. Throughout this paper we denote this
lift by $\ug_x$ (not to be confused with $\gg_x$ of Section
\ref{sec:intro}).
It should be clear that when we refer to a neighborhood of a
cuspidal geodesic $\gamma_x$ we  mean the set of geodesics
$\gamma_y$ with $y$ in a neighborhood of $x$.  By an \em interval
\em about $\ug_x$ we mean the set of geodesics $\{\ug_y : y\in
(a,b)\}$, for some $a < x < b$.

\begin{defn}[Horodisk]

Let $p \in \sinfty$.  If $p = \infty$ then an \em open horodisk
\em centered at $p$ is a set of the form $\{z| \mbox{Im}{z}
> h\}$ (for some fixed $h>0$).  For $p \in\mathbb R$ an open
horodisk centered at $p$ is a Euclidean open disk of Euclidean
radius $r$ centered at $p+ir$ (for some $r >0$). A closed horodisk
is obtained by taking the closure in $\ucover\surf$ of an open
horodisk.  Thus a closed horodisk is a closed upper half plane or
a closed disk with a point removed from its boundary.
\end{defn}

The fundamental group $\pi_1(\surf)$ can be viewed naturally as a
discrete, fixed-point free subgroup of the orientation preserving
isometries of $\ucover\surf$. The following classification is
well known:

\begin{lem}[Classification of hyperbolic geodesics]
\label{lem:classification of geodesics}

Let $f:\ucover\surf \to \ucover\surf$ be an orientation preserving
isometry, other than the identity. Then exactly one of the
following holds:

\begin{enumerate}
    \item (Elliptic) $f$ has a fixed point in the interior of
    $\ucover\surf$.  Such $f$ cannot be an element of $\pi_1(\surf)$.

    \item (Parabolic) $f$ has exactly one fixed point and this
    fixed point is on $\sinfty$.  Locally the
    fixed point is an attractor on one side and repeller on the
    other.  If $f \in \pi_1(\surf)$ then the
    fixed point is said to be a lift of the cusp and a
    sufficiently small horodisk centered at the fixed point is an
    infinite cyclic cover of a neighborhood of the cusp.

    \item (Hyperbolic) $f$ has exactly two fixed points and these
    fixed points are on $\sinfty$.  One fixed point is an attractor
    and the other repeller.  The geodesic connecting the fixed
    points is called the \em axis \em of $f$ (denoted in this paper
    $\axis{f}$) and $f$ acts on $\axis{f}$ by translation.
\end{enumerate}
\end{lem}

Axes of hyperbolic isometries correspond to closed geodesics on
$\surf$: a hyperbolic isometry $f \in \pi_1(\surf)$ generates a
cyclic group that acts on the axis $\axis{f}$ and the quotient is
a (not necessarily simple) closed geodesic on $\surf$. Conversely,
every closed geodesic on $\surf$ has a lift to $\mathbb  H^2$ that
is an (open) geodesic, say $\axis{f}$. Then $\mathbb Z$ (here, the
fundamental group of the circle) acts on $\axis{f}$ to give the
closed geodesic.  Any non-trivial element of this cyclic group is
a hyperbolic isometry whose axis is $\axis{f}$.

\begin{lem}[Being a  lift of a cusp is well-defined]
No point on $\sinfty$ is a fixed point of both a hyperbolic and a
parabolic transformation in $\pi_1(\surf)$.
\end{lem}

\begin{proof}
Let $p$ be a fixed point of a parabolic isometry.  As remarked
above a small horodisk about $p$ is an infinite cyclic cover of a
neighborhood of the cusp.  Therefore any geodesic that has $p$ as
one of its endpoints projects to a non-compact geodesic on
$\surf$.  On the other hand let $q$ be a fixed point of a
hyperbolic isometry $f$; then $q$ is the endpoint of the axis of
$f$ which projects to a closed curve on $\surf$.
\end{proof}

Since $\surf$ has finite area, fundamental domains in
$\ucover\surf$ must become arbitrarily small (in the Euclidean
metric) as they approach $\Rr$; we therefore have:

\begin{lem}
\label{lem:lifts of c are dense} The lifts of the cusp are dense
in $\sinfty$.
\end{lem}

Note too that since the set of lifts of the cusp is countable, the
complement of this set is also dense.

Like every study of simple geodesics, we need:

\begin{lem}\label{lem:simplicity closed}
Simplicity is a closed condition.
\end{lem}

\begin{proof}

Let $\ucover\gamma$ be a non-simple geodesic.  Then
$\ucover\gamma$ intersects one of its images, say
$\ucover\gamma'$.  A small perturbation of $\ucover\gamma$ gives a
small perturbation of $\ucover\gamma'$ and so the two geodesics
still intersect. Therefore non-simplicity  is an open condition
and simplicity is a closed condition.
\end{proof}

We end this section with a well-known lemma that relates geodesics
to topological properties of curves.  A curve on $\ucover\surf$ is
called proper if its ends are at infinity.  Two proper curves are
said to be \em properly homotopic \em if there is a homotopy
taking one to the other such that the ends remain at infinity at
all times. All homotopies considered in this paper are proper.

\begin{lem}
\label{lem: geodesic homotpis to simple is simple}

Let $\ucover\alpha$ be a simple curve that is properly homotopic
to a geodesic $\ug$.  The $\ug$ is simple as well.
\end{lem}

\begin{proof}
If $\ug$ is not simple then there exists a $g \in \pi_1(\surf)$ so
that $\ug \cap g(\ug) \neq \emptyset$; equivalently,  the
endpoints of $\ug$ separate the endpoints of $g(\ug)$.  Since the
homotopy is proper, endpoints do not move.  Thus $\ucover\alpha$ must
intersect $g(\ucover\alpha)$ and $\ucover\alpha$ is not simple.
\end{proof}

\section{The Proof of McShane's Theorem}
\label{sec:mainProof}

We pause for the proof of McShane's theorem, postponing several
required technical lemmas.

\begin{proof}[Proof of McShane's Theorem]

In Section~\ref{sec:deadzone} we show that every simple bicuspidal
geodesic $\gg$ lies in an open interval of otherwise non-simple
geodesics (\ref{thm:deadzone}); this interval is called the \em
deadzone \em of $\gg$ and $\gg$ is called the \em center \em of
the deadzone. Next we show that the endpoints of each deadzone are
simple (\ref{pro:deadzone and simplicity}) and unicuspidal
(\ref{prop:ends of deadzone are unicuspidal}). Consequently, the
deadzones are disjoint.

In Section~\ref{sec:non simple in DZ} we show that every
non-simple geodesic lies in a deadzone (\ref{prop:non-simple in
deadzone}).  In particular we have that any bicuspidal geodesic
lies in a deadzone (a non simple geodesic lies in a deadzone
by~\ref{prop:non-simple in deadzone} and a simple geodesic is the
center of a deadzone); hence since the lifts of the cusp are dense
(\ref{lem:lifts of c are dense}), so too is the union of the
deadzones. Letting $K\subset \Ss^1$ parameterize the set of simple
unicuspidal geodesics  we have that $K$ is exactly the complement
of the union of the deadzones. This has two immediate
consequences: first every interval of $\Ss^1 \setminus K$ contains
exactly one point of $E$ and these points correspond to the
bicuspidal geodesics (the centers of the deadzones). Second, $K$
is totally disconnected and closed.  In order to show $K$ is a
Cantor set we need only show $K$ is perfect---that every point
$x \in K$ is a limit point of $K \setminus \{x\}$.

In Section~\ref{sec:convergence onto rational} we show that the
endpoints of $K$ are limit points (\ref{prop:endpoints of
deadzones are limit points}). It follows that the remaining points
of $K$ are limit points as well: choose $x\in K\subset\Ss^1$ that
is not an endpoint; since the lifts of the cusp are dense
(\ref{lem:lifts of c are dense}) in any open interval $(a,b)$
containing $x$ there exist a lift of the cusp, say $y$. Since
$\ug_y$ is bicuspidal it is contained in a deadzone, say $V$.  By
definition $K\cap V=\emptyset$ (in particular $x$ is not$ V$), and
by assumption $x$ is not an endpoint; therefore an endpoint of $V$
must also lie in $(a,b)$ (in fact between $x$ and $y$). Thus $x$
is a limit point of $K$ and we see that $K$ is perfect as
required.

In Section~\ref{sec:classification} we show that the simple
cuspidal geodesics on $\surf$ that spiral into simple closed
geodesics are exactly the projections of endpoints of deadzones
(\ref{prop:denpoint project to spiraling geodesics}); with this we
complete the proof of McShane's Theorem.
\end{proof}

\section{Deadzones}
\label{sec:deadzone}

We now show that any simple bicuspidal geodesic $\gamma$ has a
neighborhood that contains no simple geodesic except
$\gamma$ itself.  Before stating the theorem we describe the
geometry of a simple bicuspidal oriented geodesic $\gamma$.   Let
$\ucover\gamma_x$ be a vertical lift of $\gamma$ oriented
downwards and ending at the point $x \in \Rr$.    A horodisk
centered at $x$ is an infinite cyclic cover of the cusp. In it
there are infinitely many preimages of $\gamma$, with their
orientations alternating in and out. Therefore, the two preimages
next to $\ucover\gamma_x$ are both oriented away from $x$;
concentrating on the right, we denote the preimage of $\gamma$ on
the right by $\ucover\gamma_1$; the left side is treated
similarly.  Let $f \in \pi_1(\surf)$ be the isometry taking
$\ucover\gamma_0 := \ucover\gamma_x$ to $\ucover\gamma_1$. For
$i=2,3,4,\ldots$ let $\ucover\gamma_i = f^i(\ucover\gamma_0)$ and
let $x_i$ be the terminal point of $\ucover\gamma_i$ (thus $x_0 =
x$ and $x_{i+1} = f(x_i)$). Let $\rtp = \lim_{i \to \infty} x_i$;
obviously, $f(\rtp) = \rtp$.  Similarly, by considering the left
side we construct $\ltp$. The geodesics $\ucover\gamma_i =
f^i(\gg_0)$ inherit an orientation from $\gg_0$.  We now show that
the geodesics $\ucover\gamma_i$ look just right, \ie\ just look
right:

\begin{figure}[h]
\centerline{  \includegraphics[]{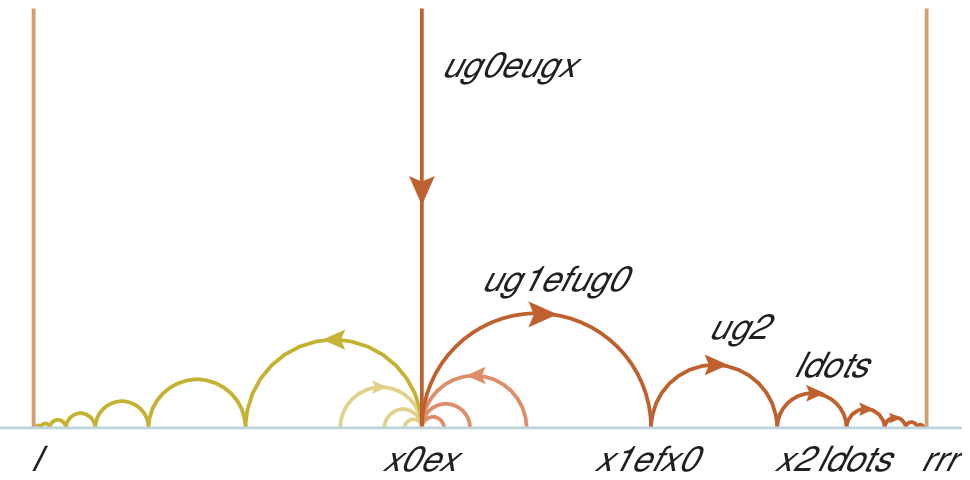}}
    \caption{A vertical lift $\ug_x$ and the corresponding chain of
    $\{\ug_i\}$ to the right of $\ug_x$,  the map $f$, and the
    interval $(\ltp,\rtp)$.}\label{fig:nomenclature}
\end{figure}

\begin{prop}
\label{pro:maximal chain}

Let $\ug_x$ be a vertical simple bicuspidal geodesic.  For the geodesics
$\ucover\gamma_i$ and the points $x_i$ constructed above we have:

\begin{enumerate}
    \item No image of $\ucover\gamma$ starts or ends at $x_i$ and
    lies between $\ucover\gamma_i$ and $\ucover\gamma_{i+1}$,
    ($i=0,1,2,\dots$).
    \item The geodesics $\ucover\gamma_i$ are all oriented consistently
    to the right (for $i \geq 1$).
    \item For any $i \neq j$ $\ucover\gamma_i \cap \ucover\gamma_j =
    \emptyset$.
\end{enumerate}
\end{prop}

\begin{proof}
\begin{enumerate}
    \item If there exists a geodesic at $x_i$ that lies between
    $\ucover\gamma_i$ and $\ucover\gamma_{i+1}$, then its image under
    $f^{-i}$ is a geodesic between $\ucover\gamma_0$ and $\ucover
    \gamma_1$, contradicting our choice of $\ucover\gg_1$.
    \item For $\ug_1$ this is true by construction; if $i > 1$ and
    $\ug_i$ is oriented to the left, it will have to either
    intersect $\ug_j$ for some $j<i$ or terminate at $x_j$ (for
    some $j<i$).  The former is impossible since $\ug_x$ is
    embedded and the latter contradicts case (1) above. Moreover
    $\ug_i$ cannot be vertical, for if it were, $f^{i+1}(\infty)=\infty$
    for $i > 1$, while $f(\infty) \neq \infty$; but then
    $f$ would be elliptic.
    \item Trivial since $\ug_x$ is embedded.
\end{enumerate}
\end{proof}

\noindent{\bf Remark.} Proposition~\ref{pro:maximal chain} holds
when $\ug_x$ is not simple, though we do not show or require this
in this paper.

Every simple bicuspidal geodesic has a neighborhood of non-simple
geodesics:

\begin{thm}\label{thm:deadzone}

Let $\ucover\gg_x$ be a simple bicuspidal geodesic.  If $c \in
(\ltp,x) \cup (x,\rtp)$ then $\ucover\gg_c$ is not simple.
\end{thm}

\begin{figure}
 \centerline{\includegraphics[]{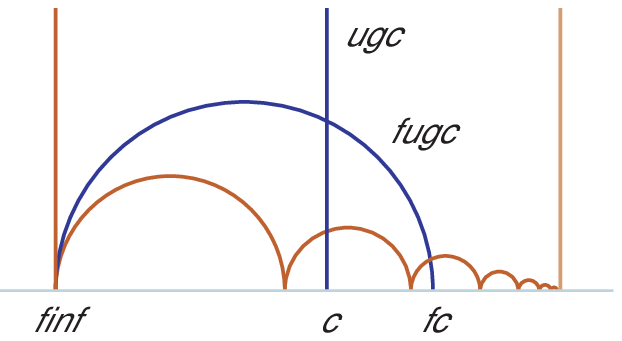}}
  \caption{No geodesic $\ug_c$, $c\in(\ltp,x)\cup(x,\rtp)$ is simple.}
  \label{fig:deadzone}
\end{figure}

\begin{proof}
Consider $\ug_c$,  $c\in(x_0,r)$ (the other case is treated
similarly). There exists an $i$ so that $c\in(x_{i-1},x_i]$. Now
consider $f(\ug_c)$. But as illustrated in Figure~\ref{fig:deadzone} the endpoints of $\ug_c$ are $\infty$ and
$c$ and the endpoints of $f(\ug_c)$ are $x_0$ and
$f(c)\in(x_i,x_{i+1}]$. These two curves must cross and $\ug_c$
cannot be simple.
\end{proof}

\begin{defn}[Deadzone]\label{def:the deadzone}
For a simple bicuspidal geodesic $\ug_x$, we call the interval
$(\ltp,\rtp)$ the \ddefine{deadzone} of $\ug_x$ and call $\ug_x$
the \ddefine{center} of its deadzone.
\end{defn}

We next  prove:

\begin{prop}\label{prop:ends of deadzone are unicuspidal}
If $\ug_x$ is bicuspidal and simple, then neither $\rtp$  nor
$\ltp$ is a lift of the cusp (and so  $\ug_{\rtp}, \ug_{\ltp}$ are
unicuspidal).
\end{prop}

\begin{figure}[h]
\centerline{\includegraphics[]{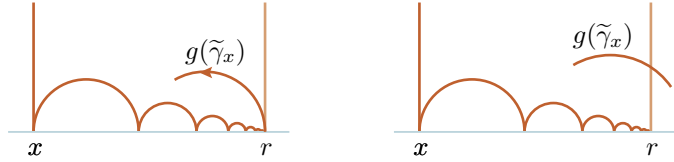} }
  \caption{Left: Proposition \ref{prop:ends of deadzone are unicuspidal}; right: Lemma \ref{lem:simple means no intersection}}
  \label{fig:figure 44}
\end{figure}

\begin{proof}
(For this proposition and the next lemma, see Figure \ref{fig:figure 44}.)
Suppose for a contradiction that $\rtp$ is a lift of the cusp;
$\ltp$ is treated similarly.
Since a horodisk about $\rtp$ is an infinite cyclic cover of the
cusp there exists an image $g(\ug_{x})$, $g\in\pi_1(\surf)$
starting at $\rtp$ oriented away from the cusp to the left.  By
(1) of Proposition~\ref{pro:maximal chain} $g(\ug_{x})$ cannot
start or terminate at any of the points $x_i$.  It must therefore
intersect $\ucover\gg_i$ (for some $i$) contradicting the
assumption that $\ucover\gamma_x$ is simple.
\end{proof}

We next show that the endpoints of a deadzone are simple:

\begin{lem} \label{lem:simple means no intersection}
Let $\ug_x$ be a simple biscuspidal geodesic. No image of $\ug_x$
meets $\ug_{\rtp}, \ug_{\ltp}$.
\end{lem}

\begin{proof}
Suppose for contradiction that some image of $\ug_x$ crosses
$\ug_{\rtp}$ ($\ug_{\ltp}$ is treated similarly). Some endpoint of
this image must lie to the left of $\rtp$, but this endpoint
cannot be any of the $x_i$'s by (1) of
Proposition~\ref{pro:maximal chain}; on the other hand this image
of $\ug_x$ cannot cross any $\ug_i$ since $\ug_x$ is simple. We
have a contradiction: there is nowhere for the image to
end.\end{proof}

\begin{prop}
\label{pro:deadzone and simplicity}

Let $\ug_x$ be a simple biscuspidal geodesic. Then $\ug_{\rtp},
\ug_{\ltp}$ are simple.
\end{prop}

\begin{figure}
 \centerline{\includegraphics[]{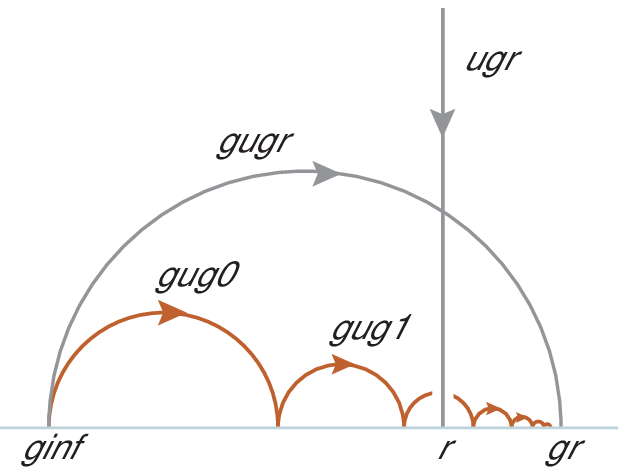}}
  \caption{$\ug_r$, $\ug_l$ are simple.}
  \label{fig:ends simple}
\end{figure}

\begin{proof}
We show $\ug_{\rtp}$ is simple; $\ug_{\ltp}$ is treated similarly.
Suppose $\ug_{\rtp}$ is non-simple. Let $g \in \pi_1(\surf)$ be an
isometry so that $g(\ug_{\rtp}) \cap\ug_{\rtp} \neq \emptyset$. If
$g(\ug_{\rtp})$ is oriented to the left then since $g$ is
orientation preserving $g^{-1}(\ug_{\rtp})$ intersects
$\ug_{\rtp}$ and will be oriented to the right. Thus without loss
of generality we may assume that $g(\ug_{\rtp})$ is oriented to
the right.  Hence $g(\infty)$ lies to the left of $\rtp$ and
$g(\rtp)$ lies to the right.  The geodesics $\{g(\ucover
\gamma_i)\}$ form a chain of geodesics from $g(\infty)$ to
$g(\rtp)$, \em beneath \em the image of $g(\ug_{\rtp})$ (see
Figure~\ref{fig:ends simple}). By Proposition~\ref{prop:ends of
deadzone are unicuspidal}, $\rtp$ is not a lift of the cusp and
therefore cannot be the endpoint of $g(\ug_i)$ for any $i$. Hence
$g(\ug_i)$ must intersect $\ug_{\rtp}$. But this contradicts
Lemma~\ref{lem:simple means no intersection}.
\end{proof}

This is now trivial:
\begin{cor}
The deadzones are disjoint.
\end{cor}

\section{Each non-simple geodesic lies in a deadzone}
\label{sec:non simple in DZ}

\begin{defn}

Let $\ug_x$ be a non-simple geodesic.   A \ddefine{point of
return} is a point $p=\ug_x\cap g(\ug_x)$ for some $g\in
\pi_1(\surf)$, such that $p$ is below $g^{-1}(p)$ (\ie\
$\Im(p)<\Im(g^{-1}(p)$). (Note that $g^{-1}(p)$ also lies on
$\ug_x$.)  The \em highest  point of return \em (if it exists) is a
point of return that has the largest imaginary value among all
points of return on $\ug_x$.
\end{defn}

\begin{lem}
\label{lem:existnce of point of return}

Each non-simple bicuspidal geodesic has a highest point of return.
\end{lem}

\begin{proof}
Trivial, since a bicuspidal geodesic meets only finitely many of
its images.
\end{proof}

\begin{rem}
\label{rem:highest point of return always exists} In fact, it is
not hard to show that any non-simple cuspidal geodesic has a
highest point of return.
\end{rem}

\begin{lem}\label{lem:ginf simple}

Let  $\ug$ be non-simple, with highest point of return
$p=\ug \cap g(\ug)$, $g\in\pi_1(\surf)$. Then
$\ug_{g(\infty)}$ is simple and bicuspidal.\end{lem}

\begin{proof}\marginpar{lables}
Since both $\infty$ and $g(\infty)$ are images of the cusp,
$\ug_{g(\infty)}$ is bicuspidal. Let $\beta_1$ be the portion of
$\ug$ above $g^{-1}(p)$ and let $\beta_2$ be the portion between
$g^{-1}(p)$ and $p$ (see Figure \ref{fig:finding simple
geodesic}).  Let $\alpha$ be the curve obtained from $\beta_1 \cup
\beta_2\cup g(\beta_1)$ after slight perturbation to avoid $p$.
Then since $p$ is the highest point of return $\alpha$ meets no
non-trivial image of itself and is therefore simple. By Lemma
\ref{lem: geodesic homotpis to simple is simple} $\ug_{g(\infty)}$
is simple as well.
\end{proof}

\begin{figure}[h]
\centerline{\includegraphics[]{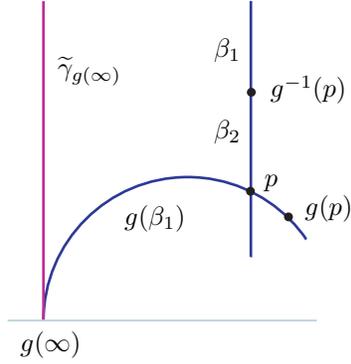} }
  \caption{A point of return $p$, the curve $\alpha=\beta_1\cup\beta_2\cup
                            g(\beta_1)$, and $\ug_{g(\infty)}$}
  \label{fig:finding simple geodesic}
\end{figure}

\begin{lem}
\label{lem:between ginf non-simple} Let  $\ug_x$ be non-simple,
with highest point of return $p=\ug_x\cap g(\ug_x)$,
$g\in\pi_1(\surf)$. Then for any $c$ between $x$ and
$g(\infty)$, the geodesic $\ug_c$ is non-simple.
\end{lem}

\begin{proof}

Let $q=\ug_c\cap g(\ug_x)$, which exists by hypothesis.
We first claim that $\Im(g^{-1}(q))>\Im(q)$: the distance from $p$ to
$q$, along $g(\ug_x)$, is the same as the distance from $g^{-1}(p)$ to
$g^{-1}(q)$ along $\ug_x$. But $\Im(g^{-1}(p))>\Im(p)$ by hypothesis, and
moreover $\ug_x$ is a vertical geodesic. Consequently,
$\Im(g^{-1}(q))>\Im(q)$ as illustrated at left in Figure~\ref{fig: fig5}

\begin{figure}[h]
{\includegraphics[]{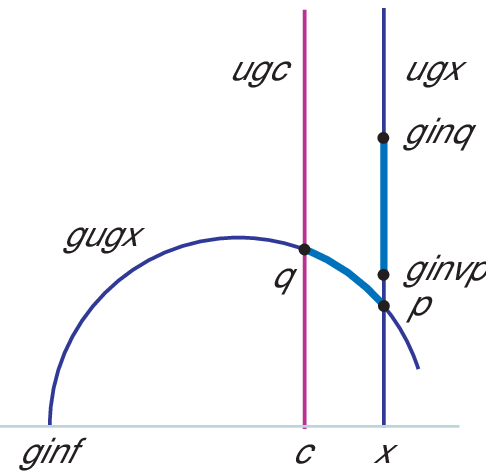} } \hfill{\includegraphics[]{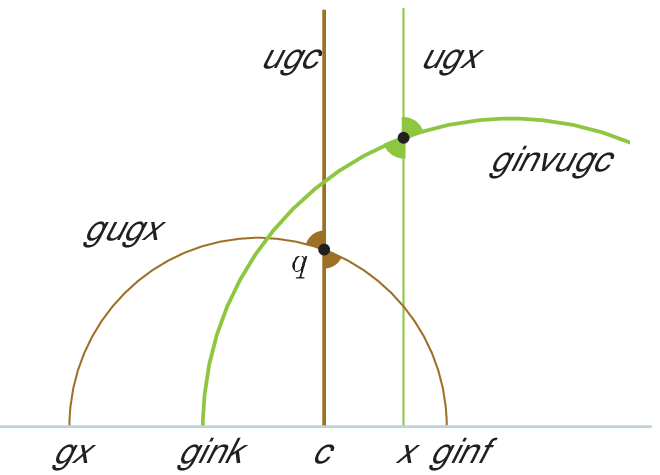}}
  \caption{At left, $\Im(g^{-1}(q))>\Im(q)$.
  At right, the similar configurations $A$ (brown) 
  and $g^{-1}(A)$
  (green). 
  }
   \label{fig: fig5}
\end{figure}

Next consider the following two
configurations:$$A=\{\ug_c,g(\ug_x), q\}$$
$$g^{-1}(A)=\{g^{-1}(\ug_c),\ug_x,g^{-1}(q)\}$$ In the upper-half plane
model of hyperbolic space, both configurations consist of a
vertical ray, a half-circle and their point of intersection.
Moreover, since $g$ is conformal, the ray and the half-circle meet
at the same angle in the two configurations. That is, the two
configurations are {\em similar} in the Euclidean sense. In fact
though, since $\Im(g^{-1}(q))>\Im(q)$, configuration $g^{-1}(A)$
is larger than configuration $A$ in the Euclidean sense. Moreover,
by examination we see that the orientations of the two
configurations are reversed.  Consequently,
$|x-g^{-1}(c)|>|g(\infty)-c|>|x-c|$ and $c$ must lie between
$g^{-1}(c)$ and $x$.
That is, $g^{-1}(\ug_c)$ must cross $\ug_c$, and so $\ug_c$ is not
simple.
\end{proof}

\begin{lem}
\label{lem:non-simple bicusp in deadzone} Each non-simple
bicuspidal geodesic $\ug$ lies in a deadzone.
\end{lem}

\begin{proof}
Let $\ucover\gg_x$ be a non-simple bicuspidal geodesic and let
$g\in\pi_1(\surf)$ be such that $p=\ug_x\cap g(\ug_x)$ is the
highest point of return; let $\ug_{g(\infty)}$ be as above. By
Lemma~\ref{lem:ginf simple}, the geodesic $\ug_{g(\infty)}$ is
simple and bicuspidal and so (as discussed in
Section~\ref{sec:deadzone}) has a deadzone (say $V$) of non-simple
geodesics, bounded by a pair of simple geodesics
(\ref{pro:deadzone and simplicity}).  By Lemma~\ref{lem:between
ginf non-simple} no endpoint of $V$ lies between $\ug_{g(\infty)}$
and $\ug_x$; by assumption $\ug_x$ is bicuspidal and is not an
endpoint; and so $\ug_x$ must lie in $V$.
\end{proof}

\noindent{\bf Remarks.} (1) Our proof is constructive, finding
$\ug_{g(\infty)}$ the
    center of the deadzone.  Compare with \cite[Lemma 5.8]{tan-main}.

    (2) Note that we proved that a non-simple cuspidal geodesic that has a
    highest point of return lies in a deadzone. By Remark
    \ref{rem:highest point of return always exists} this proof holds
    for all non-simple cuspidal geodesics. Perhaps the proof below is
    slightly easier.

\begin{prop}
\label{prop:non-simple in deadzone}

Every non-simple geodesic lies in a deadzone.
\end{prop}

\begin{proof}   Let $\ug_x$ be any non-simple geodesic; by Lemma
\ref{lem:non-simple bicusp in deadzone} we may assume $\ug_x$ is
unicuspidal.  By Lemma~\ref{lem:simplicity closed} there exists an
open interval $(a,b) \subset \Rr$ containing $x$  such that for
any $c\in (a,b)$, $\ug_c$ is non-simple as well. Since lifts of
the cusp are dense, there exists $c\in (a,b)$ such that $\ug_c$ is
bicuspidal and non-simple. By  Lemma~\ref{lem:non-simple bicusp in
deadzone} $\ug_c$ is in the deadzone of some bicuspidal simple
geodesic; the endpoints of this deadzone are simple
(\ref{pro:deadzone and simplicity}) and hence cannot lie in
$(a,b)$ and so in fact all of $(a,b)$ (including $x$) lies in this
deadzone as well.
\end{proof}

\section{Endpoints of deadzones are the limit of simple geodesics }
\label{sec:convergence onto rational}

We now show that there are simple bicuspidal geodesics arbitrarily
close to every  endpoint of each deadzones.  In the following
proposition we consider a geodesic $\gamma_\rtp$ at the right end
of a deadzone; the proof for a geodesic at the left end is
precisely the same.

\begin{prop}
\label{prop:endpoints of deadzones are limit points}

Let $\ug_\rtp$ be the right endpoint of a deadzone. Then for every
$\epsilon>0$, there exists an $x\in (\rtp,\rtp+\epsilon)$ so that
$\ug_{x}$ is simple and bicuspidal.
\end{prop}

\begin{proof}

Let $f\in\pi_1(\surf)$ be a transformation fixing $\rtp$; by
Proposition~\ref{prop:ends of deadzone are unicuspidal} $f$ is
hyperbolic.  By replacing $f$ with $f^{-1}$ is necessary we may
assume $r$ is an attractor.  Because $\ug_\rtp$ is the {\em right}
endpoint of a deadzone, the axis $\axis{f}$ of $f$ has its other
endpoint to the right of $\rtp$ and the images of $\ug_\rtp$ under
$f^n$, $n\in\Zz$ are as illustrated in Figure~\ref{fig:fig6a}.

\begin{figure}[h]
\centerline{\includegraphics[]{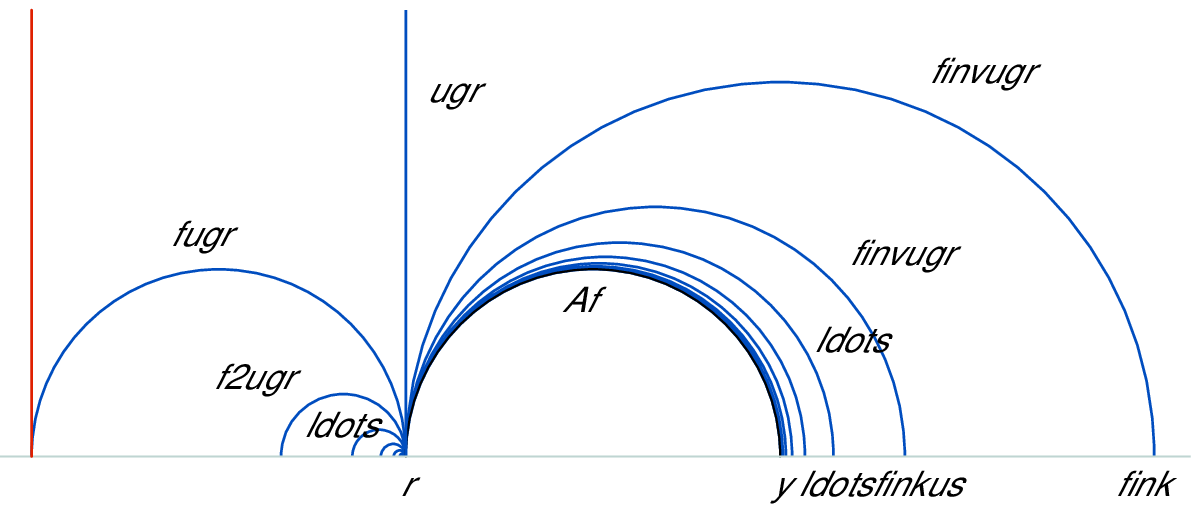} }
  \caption{$\ug_{\rtp}$ and its images under $f^n$, the axis
  $\axis{f}$ and fixed point $y$}
   \label{fig:fig6a}
\end{figure}

The point $r $ is an attracting fixed point of $f$. Denote the
other fixed point of $f$ by $y$.  Then $f^{-n}(\infty)$ limits
onto $y$.  Thus we see that $f^{-n}(\ug_{\rtp})$ limits onto
$\axis{f}$; by Lemma \ref{lem:simplicity closed} $\aaxis_{f}$ is
consequently simple.

\begin{figure}[h]
{\includegraphics[]{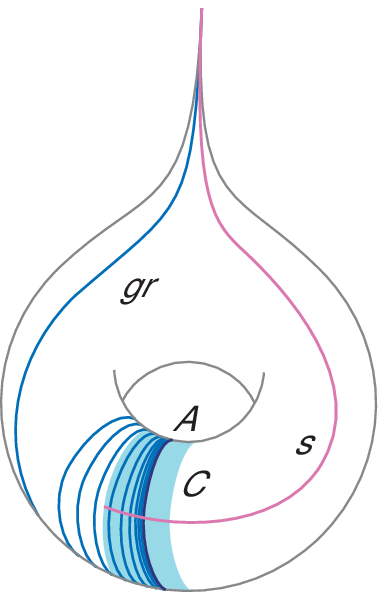} } \hfill{\includegraphics[]{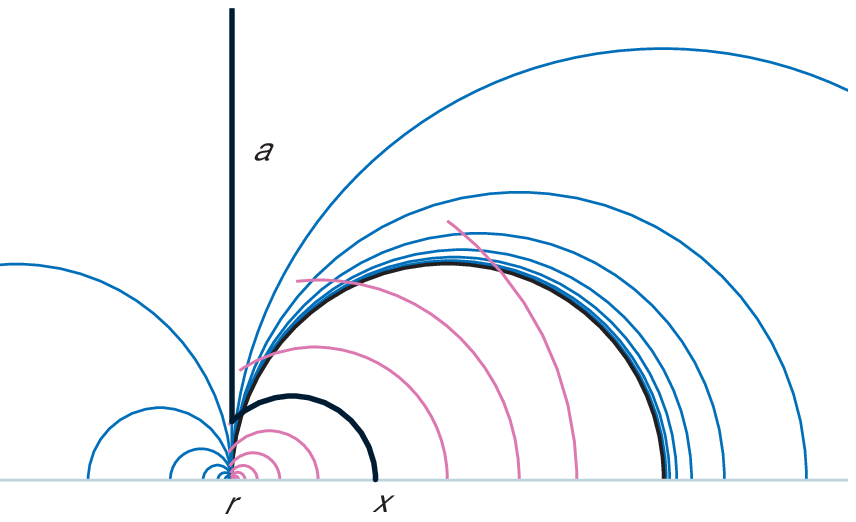}}
\caption{At left, the curve $\sigma$, crossing the neighborhood
$C$ of $\aaxis$ on $\surf$. At right the construction of $\alpha$
from  $\ug_r$ and a suitable lift of $\sigma$.}
   \label{fig:fig6b}
\end{figure}

Let $\aaxis$ be the image of $\axis{f}$ on $\surf$, $\gg_r$ the
image of $\ucover\gg_r$ and $C$ a neighborhood of $\aaxis$ as
illustrated at left in Figure~\ref{fig:fig6b}. Take $\sigma$ to be
any embedded curve that leaves the cusp, misses
$\gamma_{\rtp}\setminus C$, and crosses $C$, meeting $\aaxis$
exactly once, entering on the opposite side from $\gamma_{\rtp}$
as pictured.

Choose any lift $\ucover{\sigma}$ of $\sigma$ that meets
$\axis{f}$, and let $x'$ be the image of the cusp at the end of
$\ucover{\sigma}$.  Note that $x'$ lies beneath $\axis{f}$ to the
right of $\rtp$, since $\sigma$ approaches $\aaxis$ on the
opposite side from $\gamma$.  Recalling that $\rtp$ is an
attracting point for $f$ and that $\sigma$ meets $\gamma$
infinitely often, we note that for any sufficiently large
$n\in\Nn$, $f^n(\ucover\sigma)\cap \ug_\rtp\neq\emptyset$ and
$f^n(x')\in(\rtp,\rtp+\epsilon)$. Choose any such $n$ and let
$x=f^n(x')$.

We now claim that the  vertical geodesic $\ug_{x}$ is simple,
completing the proof of the lemma:  Let $\ucover{\alpha}$ consist
of the portion of $f^n(\ucover{\sigma})$ between $x$ and
$\ug_{\rtp}$ and the portion of $\ug_{\rtp}$ above
$f^n(\ucover{\sigma}) \cap \ug_r$. Then $\ucover{\alpha}$ is
simple and homotopic to $\ug_x$; consequently by Lemma \ref{lem:
geodesic homotpis to simple is simple} $\ug_x$ is simple as well.
\cmp{We should add this lemma; it is also useful for Lemma 5.3}
\end{proof}

\section{Endpoint of deadzones spiral onto simple closed curves}
\label{sec:classification}

We finished proving that $K$ is a Cantor set and $\Ss^1 \setminus
K$ has exactly one point in each complimentary interval; those
points are exactly the points corresponding to bicuspidal
geodesics.  To complete McShane's Theorem we need to understand
the endpoints of the deadzones.  We say that an oriented geodesic
$\gamma$ on $\surf$ \em\ spirals onto \em\ a (not necessarily
simple) closed curve $\alpha$ if for every $\epsilon > 0$ there
exists $t_0$ so that for any $t > t_0,$ $\gamma(t)$ is $\epsilon$
close to $\alpha$.

\begin{prop}
\label{prop:denpoint project to spiraling geodesics}

A simple cuspidal geodesic $\ug_x$ projects to a geodesic that
spirals onto a simple closed geodesic if and only if $\ug_x$ is an
endpoint of a deadzone.
\end{prop}

\begin{proof}
\noindent{\bf Claim:} $\ug_x$ projects to a geodesic that spirals
onto a (not necessarily simple) closed curve if and only if $x$ is
a fixed point of a hyperbolic isometry.

To prove the claim, first assume that $x$ is the fixed point of a
hyperbolic isometry (say $f$).  Then $x$ is an endpoint of the
axis of $f$, $\aaxis_f$.  Since $\ug_x$ approaches $\aaxis_f$
asymptotically and $\aaxis_f$ projects to a closed geodesic, the
claim holds in this case.  Conversely, suppose that the projection
of $\ug_x$ spirals onto a closed geodesic, say $\alpha$.  Then
$\ug_x$ gets arbitrarily close to some lift of $\alpha$, say
$\ucover\alpha$. Then $\ucover\alpha$ is an axis of a hyperbolic
isometry establishing the claim.

If $\ug_x$ is the endpoint of a deadzone then the projection of
$\ug_x$ spirals onto the projection of $\aaxis_{f}$.  As we saw in
Section~\ref{sec:convergence onto rational} $\ug_x$ accumulates
onto $\aaxis_f$ and hence $\aaxis_f$ is simple
(\ref{lem:simplicity closed}).  Thus the projection of $\gg_x$
spirals onto the projection of $\aaxis_f$, establishing the
proposition in this case.

Conversely, suppose that $\ug_x$ is a simple geodesic that
projects to a geodesic that spirals onto a simple closed geodesic.
Then some lift of the closed geodesic ends at $x$; as before, we
denote the other endpoint of this lift by $y$ and note that since
$y$ is not a lift of the cusp $y \neq \infty$.  Without loss of
generality we assume that $y$ is to the right of $x$ (see
Figure~\ref{fig:7a}).
\begin{figure}[h]
{\includegraphics[]{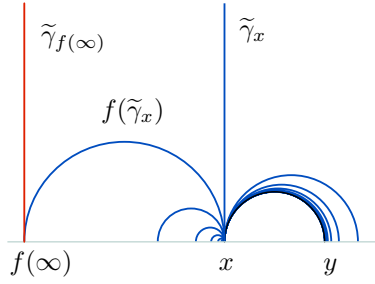}}
  \caption{Simple spiraling geodesic}
  \label{fig:7a}
\end{figure}
It is a basic exercise in planar hyperbolic geometry that if two
hyperbolic isometries in $\pi_1(\surf)$ share a fixed point then
they share the other fixed point as well.\footnote{Hint: assume
for contradiction there exist $f_1$, $f_2 \in \pi_1(\surf)$
hyperbolic isometries with axes $\aaxis_{f_1} \neq \aaxis_{f_2}$
with one fixed point in common, and show that for any $\epsilon >
0$ the projection of $\aaxis_{f_1}$ is in the
$\epsilon$-neighborhood of $\aaxis_{f_2}$ and {\it vice versa};
use this for the contradiction.}  Therefore the collection of
hyperbolic isometries that fix $x$ also fix $y$ (and hence the
geodesic connecting $x$ to $y$).  Given a collection $\{f_i \in
\pi_1(\surf)\}$ of isometries that fix $x$ and $y$, for any $p \in
\Ss^1$ the only accumulation points of $f_i(p)$ are $x$ and $y$;
in particular, $\infty$ is not an accumulation point and there is
no accumulation point in $(-\infty,x)$. Therefore there exists a
unique $f \in \pi_1(\surf)$ so that for any $f \neq f' \in
\pi_1(\surf)$ with $f'(x) = x$, $f(\infty) <
f'(\infty)$.\marginpar{FIX}

\begin{figure}[h]
{\includegraphics[]{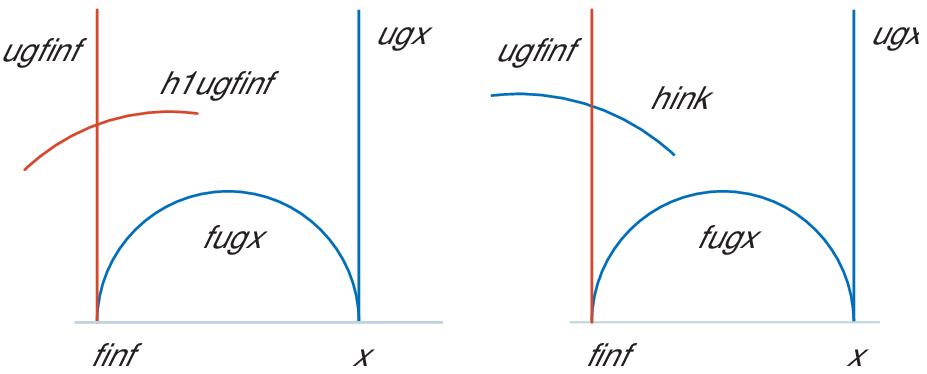}}
    \caption{$\ug_{f(\infty)}$ is  simple. }
    \label{fig:7b}
\end{figure}

Next we show that $\ug_{f(\infty)}$ is a simple geodesic; assume
for contradiction it is not.  Then for some $h_1 \in
\pi_1(\surf)$, $\ug_{f(\infty)} \cap h_1(\ug_{f(\infty)}) \neq
\emptyset$. Clearly $h_1(\ug_{f(\infty}))$ is bicuspidal, and
therefore cannot start or end at $x$.  We conclude that
$h_1(\ug_{f(\infty)})$ must intersect either $\ug_x$ of
$f(\ug_x)$.  Therefore some image of $\ug_x$ must intersect
$\ug_{f(\infty)}$, say $h_2(\ug_{x})$. Since $\ug_x$ is
simple, $h_2(\ug_{x})$ cannot intersect either $\ug_x$ or
$f(\ug_x)$.  Therefore $h_2(\ug_x)$ terminates at $x$ (it cannot
start at $x$ since $x$ is not a lift of the cusp).  Thus $h_2(x) =
x$.  But $h_2(\infty)$ lies to the left of $f(\infty)$,
contradicting our choice of $f$.
%
%

Finally,
the chain $\{f^n(\ug_{f(\infty)})\}_{n=1,2,\ldots}$ lies within
the deadzone of $\ug_{f(\infty)}$ and has right endpoint $x$.
Since $\ug_x$ is simple, we conclude that $\ug_x$ is the
right endpoint of the deadzone of  $\ug_{f(\infty)}$.

\end{proof}

\nocite{*}
\ifx\undefined\bysame
\newcommand{\bysame}{\leavevmode\hbox to3em{\hrulefill}\,}
\fi

\end{document}